\documentclass[a4paper,12pt]{article}
\usepackage{amsmath}
\usepackage{latexsym}
\usepackage{url}
\usepackage{color}
\numberwithin{equation}{section}

\def\R{{\bf R}}

\def\N{{\bf N}}

\def\d{\displaystyle}
\def\e{{\varepsilon}}

\def\p{\partial}
\def\v#1{\mbox{\boldmath $#1$}}

\newtheorem{thm}{Theorem}[section]

\newtheorem{prop}{Proposition}[section]

\title{Semilinear wave equations of derivative type with spatial weights
in one space dimension}
\author{
Shunsuke Kitamura
\footnote{
Doctor course, Mathematical Institute,
Tohoku University,
Aoba, Sendai 980-8578, Japan.
e-mail: shunsuke.kitamura.s8@dc.tohoku.ac.jp.},
Katsuaki Morisawa
\footnote{
Musashi High School and Junior High School,
1-26-1 Toyotamakami, Nerima, Tokyo, 176-8535.
e-mail: morisawa.katsuaki@musashi.ed.jp.},
\\
Hiroyuki Takamura
\footnote{Mathematical Institute,
Tohoku University,
Aoba, Sendai 980-8578, Japan.
e-mail: hiroyuki.takamura.a1@tohoku.ac.jp.}
}
\date{
\[
\begin{array}{ll}
\mbox{\footnotesize{\bf Keywords:}}
& \mbox{\footnotesize semilinear wave equation, one dimension, classical solution, lifespan}\\
\mbox{\footnotesize{\bf MSC2020:}}
& \mbox{\footnotesize primary 35L71, secondary 35B44}\\
\end{array}
\]
}
\begin{document}
\maketitle
\begin{abstract}
This paper is devoted to the initial value problems for semilinear wave equations
of derivative type with spatial weights in one space dimension.
The lifespan estimates of classical solutions
are quite different from those for nonlinearity of unknown function itself
as the global-in-time existence can be established by spatial decay.
\end{abstract}

%%%%%%%%%%%%%%%%%%%%%%%%%%%%%%%%%%%%%%%%%%%%%%%%%%
%%%%%%%%%%%%%%%%%%%%% SECTION1 %%%%%%%%%%%%%%%%%%%%%%%
%%%%%%%%%%%%%%%%%%%%%%%%%%%%%%%%%%%%%%%%%%%%%%%%%%

\section{Introduction}

\par 
In this paper, we consider the initial value problems;
\begin{equation}
\label{IVPderivative}
\left\{
\begin{array}{ll}
	\d u_{tt}-u_{xx}=\frac{|u_t|^p}{(1+x^2)^{(1+a)/2}}
	&\mbox{in}\quad \R\times(0,\infty),\\
	u(x,0)=\e f(x),\ u_t(x,0)=\e g(x),
	& x\in\R,
\end{array}
\right.
\end{equation}
where $p>1$, $a\in\R$, $f$ and $g$ are given smooth functions of compact support
and a parameter $\e>0$ is \lq\lq small enough".
We are interested in the estimate of the lifespan $T(\e)$,
the maximal existence time,
of classical solutions of (\ref{IVPderivative}).
Our result is the following;
\begin{equation}
\label{lifespan}
\begin{array}{l}
T(\e)\sim
\left\{
\begin{array}{ll}
C\e^{-(p-1)/(-a)} & \mbox{for}\ a<0,\\
\exp\left(C\e^{-(p-1)}\right) & \mbox{for}\ a=0,
\end{array}
\right.
\\
T(\e)=\infty\qquad\mbox{for}\quad a>0.
\end{array}
\end{equation}
Here we denote the fact that there are positive constants,
$C_1$ and $C_2$, independent of $\e$ satisfying $A(\e,C_1)\le T(\e)\le A(\e,C_2)$
by $T(\e)\sim A(\e,C)$.
We note that (\ref{lifespan}) is established for classical solutions when $p\ge2$,
while we have to consider $C^1$ solutions of associated integral equations to (\ref{IVPderivative})
in case of $1<p<2$.
When $a=-1$, the upper bounds in (\ref{lifespan}) are already obtained by Zhou \cite{Zhou01},
while  the lower bounds are verified only for integer $p$ by general theory
which is studied by Li, Yu and Zhou \cite{LYZ91, LYZ92}.
We see that (\ref{lifespan}) is similar to the one for the time-weighted nonlinear terms of
unknown function itself by Kato, Takamura and Wakasa \cite{KTW19}
in sense that there is a possibility to obtain the global-in-time existence in spite of one dimension.
For such an equation, the lifespan estimates are classified into two cases according to
the value of the total integral of the initial speed.
But (\ref{lifespan}) has no classification whatever it is.
This is due to the fact that Huygens' principle is always available
for the time derivative of the solution of the free wave equation.

\par
In fact, let us compare (\ref{IVPderivative}) with
\begin{equation}
\label{IVP}
\left\{
\begin{array}{ll}
	\d u_{tt}-u_{xx}=\frac{|u|^p}{(1+x^2)^{(1+a)/2}}
	&\mbox{in}\quad \R\times(0,\infty),\\
	u(x,0)=\e f(x),\ u_t(x,0)=\e g(x),
	& x\in\R.
\end{array}
\right.
\end{equation}
In our previous work \cite{KMT21}, it is established that
\begin{equation}
\label{lifespan_non-zero}
T(\e)\sim
\left\{
\begin{array}{ll}
C\e^{-(p-1)/(1-a)} & \mbox{for}\ a<0,\\
\phi^{-1}(C\e^{-(p-1)}) & \mbox{for}\ a=0,\\
C\e^{-(p-1)} & \mbox{for}\ a>0
\end{array}
\right.
\quad
\mbox{if}\ \int_\R g(x)dx\neq0,
\end{equation}
where $\phi^{-1}$ is an inverse function of $\phi$ defined by
\[
\phi(s):=s\log(2+s),
\]
and
\begin{equation}
\label{lifespan_zero}
T(\e)\sim
\left\{
\begin{array}{ll}
C\e^{-p(p-1)/(1-pa)} & \mbox{for}\ a<0,\\
\psi_p^{-1}(C\e^{-p(p-1)}) & \mbox{for}\ a=0,\\
C\e^{-p(p-1)} & \mbox{for}\ a>0
\end{array}
\right.
\quad
\mbox{if}\ \int_\R g(x)dx=0,
\end{equation}
where $\psi_p^{-1}$ is an inverse function of $\psi_p$ defined by 
\[
\psi_p(s):=s\log^p(2+s).
\]
We remark that the quantities in all the cases of (\ref{lifespan_non-zero}) are smaller than
those of (\ref{lifespan_zero}).
This work in \cite{KMT21} is an extension of the case of $a=-1$ by Zhou \cite{Zhou92},
and inspired by non-compactly supported case by Suzuki \cite{Suzuki10},
Kubo, Osaka and Yazici \cite{KOY13} and Wakasa \cite{Wakasa17}.
See Introduction in \cite{KMT21} for details.

\par
Finally, we note that our result cannot be valid with our method in this paper
if $|u_t|^p$ in (\ref{IVPderivative}) is replaced with $|u_x|^p$ even in the non-weighted case, $a=-1$,
because the blow-up part requires a positiveness of the nonlinear terms
in the level of the $x$-derivative of the unknown function.
But one can find immediately that it is impossible by the expression of $u_x$
in (\ref{nonlinear_derivative_conjugate}).
On contrast, it is possible to obtain the same lifespan estimate from below, the existence part,
along with our method.
\par
This work was almost completed when the first and second authors were in the master course of Mathematical Institute, Tohoku University
and the third author had the second affiliation with
Research Alliance Center of Mathematical Sciences, Tohoku University.
This paper is organized as follows.
In the next section, (\ref{lifespan}) is divided into two theorems,
and the preliminaries are introduced.
Section 3 is devoted to the proof of the existence part of (\ref{lifespan}).
The main strategy is the iteration method in the weighted $L^\infty$ space 
which is originally introduced by John \cite{John79}.
In Section 4, we prove a priori estimate.
Finally, we prove the blow-up part of (\ref{lifespan})
employing the method by Zhou \cite{Zhou01} in Section 5.

%%%%%%%%%%%%%%%%%%%%%%%%%%%%%%%%%%%%%%%%%%%%%%%%%%
%%%%%%%%%%%%%%%%%%%%% SECTION2 %%%%%%%%%%%%%%%%%%%%%%%
%%%%%%%%%%%%%%%%%%%%%%%%%%%%%%%%%%%%%%%%%%%%%%%%%%

\section{Preliminaries and main results}

Throughout of this paper, we assume that the initial data
$(f,g)\in C_0^2(\R)\times C^1_0(\R)$ satisfies
\begin{equation}
\label{supp_initial}
\mbox{\rm supp }f,\ \mbox{supp }g\subset\{x\in\R:|x|\le R\},\quad R\ge1.
\end{equation}
Let $u$ be a classical solution of (\ref{IVPderivative}) in the time interval $[0,T]$.
Then the support condition of the initial data, (\ref{supp_initial}), implies that
\begin{equation}
\label{support_sol}
\mbox{supp}\ u(x,t)\subset\{(x,t)\in\R\times[0,T]:|x|\le t+R\}.
\end{equation}
For example, see Appendix in John \cite{John_book} for this fact.

\par
It is well-known that $u$ satisfies the following integral equation;
\begin{equation}
\label{integral}
u(x,t)=\e u^0(x,t)+L_a(|u_t|^p)(x,t),
\end{equation}
where $u^0$ is a solution of the free wave equation with the same initial data;
\begin{equation}
\label{linear}
u^0(x,t):=\frac{1}{2}\{f(x+t)+f(x-t)\}+\frac{1}{2}\int_{x-t}^{x+t}g(y)dy,
\end{equation}
and a linear integral operator $L_a$ for a function $v=v(x,t)$ is Duhamel's term defined by
\begin{equation}
\label{nonlinear}
L_a(v)(x,t):=\frac{1}{2}\int_0^tds\int_{x-t+s}^{x+t-s}\frac{v(y,s)}{(1+y^2)^{(1+a)/2}}dy.
\end{equation}
Then, one can apply the time-derivative to (\ref{integral}) and (\ref{linear}) to obtain
\begin{equation}
\label{integral_derivative}
u_t(x,t)=\e u_t^0(x,t)+L_a'(|u_t|^p)(x,t)
\end{equation}
and
\begin{equation}
\label{linear_derivative}
u_t^0(x,t)=\frac{1}{2}\{f'(x+t)-f'(x-t)+g(x+t)+g(x-t)\},
\end{equation}
where $L_a'$ for a function $v=v(x,t)$ is defined by
\begin{equation}
\label{nonlinear_derivative}
\begin{array}{ll}
L_a'(v)(x,t):=
&\d\frac{1}{2}\int_0^t\frac{v(x+t-s,s)}{\{1+(x+t-s)^2\}^{(1+a)/2}}ds\\
&+\d\frac{1}{2}\int_0^t\frac{v(x-t+s,s)}{\{1+(x-t+s)^2\}^{(1+a)/2}}ds.
\end{array}
\end{equation}
On the other hand, applying the space-derivative to (\ref{integral}) and (\ref{linear}),
we have
\[
u_x(x,t)=\e u_x^0(x,t)+\overline{L_a'}(|u_t|^p)(x,t)
\]
and
\[
u_x^0(x,t)=\frac{1}{2}\{f'(x+t)+f'(x-t)+g(x+t)-g(x-t)\},
\]
where $\overline{L_a'}$ for a function $v=v(x,t)$ is defined by
\begin{equation}
\label{nonlinear_derivative_conjugate}
\begin{array}{ll}
\overline{L_a'}(v)(x,t):=
&\d\frac{1}{2}\int_0^t\frac{v(x+t-s,s)}{\{1+(x+t-s)^2\}^{(1+a)/2}}ds\\
&-\d\frac{1}{2}\int_0^t\frac{v(x-t+s,s)}{\{1+(x-t+s)^2\}^{(1+a)/2}}ds.
\end{array}
\end{equation}
Therefore, $u_x$ is expressed by $u_t$.
Moreover, one more space-derivative to (\ref{integral_derivative}) yields that
\begin{equation}
\label{integral_2derivative}
\begin{array}{ll}
u_{tx}(x,t)=&\d\e u_{tx}^0(x,t)\\
&\d+pL_a'(|u_t|^{p-2}u_tu_{tx})(x,t)
-(1+a)L_{a+2}'(|u_t|^px)(x,t)
\end{array}
\end{equation}
and
\begin{equation}
\label{linear_2derivative}
u_{tx}^0(x,t):=\frac{1}{2}\{f''(x+t)-f''(x-t)+g'(x+t)+g'(x-t)\}.
\end{equation}
Similarly, we have that
\begin{equation}
\label{u_{tt}}
\begin{array}{ll}
u_{tt}(x,t)=&
\d\e u_{tt}^0(x,t)+\frac{|u_t(x,t)|^p}{(1+x^2)^{(1+a)/2}}\\
&\d +p\overline{L_a'}(|u_t|^{p-2}u_tu_{tx})(x,t)
-(1+a)\overline{L_{a+2}'}(|u_t|^px)(x,t)
\end{array}
\end{equation}
and
\[
u_{tt}^0(x,t)=\frac{1}{2}\{f''(x+t)+f''(x-t)+g'(x+t)-g'(x-t)\}.
\]
Therefore, $u_{tt}$ is expressed by $u_{tx}$ and $u_t$, so is $u_{xx}$ because of
\[
\begin{array}{ll}
u_{xx}(x,t)=
&\d\e u_{xx}^0(x,t)\\
&\d+p\overline{L_a'}(|u_t|^{p-2}u_tu_{tx})(x,t)
-(1+a)\overline{L_{a+2}'}(|u_t|^px)(x,t)
\end{array}
\]
and
\[
u_{xx}^0(x,t)=u^0_{tt}(x,t).
\]

\par
First, we note the following fact.

\begin{prop}
\label{prop:continuity}
Assume that $(f,g)\in C^2(\R)\times C^1(\R)$.
Let $u_t$ be a $C^1$ solution of (\ref{integral_derivative}). 
Then,
\begin{equation}
\label{set}
w(x,t):=\int_0^tu_t(x,s)ds+\e f(x)
\end{equation}
 is a classical solution of (\ref{IVPderivative}).
\end{prop}
\par\noindent
{\bf Proof.} It is trivial that $w$ satisfies the initial condition and
\begin{equation}
\label{equality}
w_t=u_t,\quad w_{tt}=u_{tt}.
\end{equation}
Then, (\ref{integral_2derivative}) yields that
\[
\begin{array}{ll}
w_x(x,t)
&\d=\int_0^tu_{tx}(x,s)ds+\e f'(x)\\
&\d=\int_0^t\{pL_a'(|u_t|^{p-2}u_tu_{tx})(x,s)
-(1+a)L_{a+2}'(|u_t|^px)(x,s)\}ds\\
&\d\quad+\int_0^t\e u_{tx}^0(x,s)ds+\e f'(x)\\
&\d=\overline{L'_a}(|u_t|^p)(x,t)+\e u_x^0(x,t)
\end{array}
\]
because of
\[
pL_a'(|u_t|^{p-2}u_tu_{tx})(x,s)
-(1+a)L_{a+2}'(|u_t|^px)(x,s)
=\frac{\p}{\p s}\overline{L'_a}(|u_t|^p)(x,s).
\]
Therefore we obtain that
\[
\begin{array}{ll}
w_{xx}(x,t)=
&\d\e u_{xx}^0(x,t)\\
&\d+p\overline{L_a'}(|u_t|^{p-2}u_tu_{tx})(x,t)
-(1+a)\overline{L_{a+2}'}(|u_t|^px)(x,t)
\end{array}
\]
which implies, together with (\ref{u_{tt}}) and (\ref{equality}), the desired conclusion,
\[
w_{tt}-w_{xx}=\frac{|u_t|^p}{(1+x^2)^{(1+a)/2}}=\frac{|w_t|^p}{(1+x^2)^{(1+a)/2}}.
\]
\hfill$\Box$

\vskip10pt
Our result in (\ref{lifespan}) is splitted into the following two theorems.
\begin{thm}
\label{thm:lower-bound}
Assume (\ref{supp_initial}).
Then, there exists a positive constant $\e_1=\e_1(f,g,p,a,R)>0$ such that
a classical solution $u\in C^2(\R\times[0,T])$ of (\ref{IVPderivative}) for $p\ge2$,
or a solution $u_t\in C(\R\times[0,T])$ with
\[
\mbox{\rm supp}\ u_t(x,t)\subset\{(x,t)\in\R\times[0,T]:|x|\le t+R\}
\]
of associated integral equations
of (\ref{integral_derivative}) to (\ref{IVPderivative}) for $1<p<2$, exists
as far as $T$ satisfies
\begin{equation}
\label{lower-bound}
\begin{array}{l}
T\le
\left\{
\begin{array}{ll}
c\e^{-(p-1)/(-a)} & \mbox{for}\ a<0,\\
\exp\left(c\e^{-(p-1)}\right) & \mbox{for}\ a=0,\end{array}
\right.
\\
T<\infty\quad\mbox{for}\ a>0,
\end{array}
\end{equation}
where $0<\e\le\e_1$, $c$ is a positive constant independent of $\e$.
\end{thm}
\begin{thm}
\label{thm:upper-bound}
Assume (\ref{supp_initial})
and
\begin{equation}
\label{positive_non-zero}
\int_{\R}g(x)>0.
\end{equation}
Then, there exists a positive constant $\e_2=\e_2(g,p,a,R)>0$ such that
a solution $u_t\in C(\R\times[0,T])$ with
\[
\mbox{\rm supp}\ u_t(x,t)\subset\{(x,t)\in\R\times[0,T]:|x|\le t+R\}
\]
of associated integral equations (\ref{integral_derivative})
to (\ref{IVPderivative})
cannot exist whenever $T$ satisfies
\begin{equation}
\label{upper-bound}
T\ge
\left\{
\begin{array}{ll}
C\e^{-(p-1)/(-a)} & \mbox{for}\ a<0,\\
\exp\left(C\e^{-(p-1)}\right) & \mbox{for}\ a=0,
\end{array}
\right.
\end{equation}
where $0<\e\le\e_2$, $C$ is a positive constant independent of $\e$.
\end{thm}
The proofs of above theorems are given in following sections.

%%%%%%%%%%%%%%%%%%%%%%%%%%%%%%%%%%%%%%%%%%%%%%%%%%
%%%%%%%%%%%%%%%%%%%%% SECTION3 %%%%%%%%%%%%%%%%%%%%%%%
%%%%%%%%%%%%%%%%%%%%%%%%%%%%%%%%%%%%%%%%%%%%%%%%%%

\section{Proof of Theorem \ref{thm:lower-bound}}
\par
According to the observations in the previous section,
we shall construct a $C^1$ solution of (\ref{integral_derivative}) when $p\ge2$ and
a continuous solution of (\ref{integral_derivative}) when $1<p<2$.

\par
First, we shall construct a $C^1$ solution for $p\ge2$.
Let $\{U_j(x,t)\}_{j\in\N}$ be a sequence of $C^1(\R\times[0,T])$ defined by
\begin{equation}
\label{U_j}
U_{j+1}=\e u_t^0+L'_a(|U_j|^p),\ U_1=\e u_t^0.
\end{equation}
Then, in view of (\ref{integral_2derivative}), $(U_j)_x$ has to satisfy
\begin{equation}
\label{U_j_x}
\left\{
\begin{array}{l}
(U_{j+1})_x=\e u^0_{tx}+pL_a'(|U_j|^{p-2}U_j(U_j)_x)-(1+a)L'_{a+2}(|U_j|^px),\\
(U_1)_x=\e u^0_{tx},
\end{array}
\right.
\end{equation}
so that the function space in which $\{U_j\}$ will converge is
\[
X:=\{U\in C^1(\R\times[0,T]) :\ \mbox{supp}\ U\subset\{|x|\le t+R\}\},
\]
equipping the norm
\[
\|U\|_X:=\|U\|+\|U_x\|,\quad \|U\|:=\sup_{\R\times[0,T]}|U(x,t)|.
\]
\par
First we note that supp $U_j\subset\{(x,t)\in\R\times[0,T]\ :\ |x|\le t+R\}$ implies supp
$U_{j+1}\subset\{(x,t)\in\R\times[0,T]\ :\ |x|\le t+R\}$.
It is easy to check this fact by assumption on the initial data (\ref{supp_initial})
and the definitions of $L'_a$ in (\ref{integral_derivative}).

\begin{prop}
\label{prop:apriori}
Let $U\in C(\R\times[0,T])$ and {\rm supp}\ $U\subset\{(x,t)\in\R\times[0,T]:|x|\le t+R\}$. Then there exists a positive constant $C$ independent of $T$ and $\e$ such that
\begin{equation}
\label{apriori}
\|L'_a(|U|^p)\|\le CE_a(T)\|U\|^p,
\end{equation}
where
\begin{equation}
\label{E}
E_a(T):=
\left\{
\begin{array}{ll}
(T+2R)^{-a} & \mbox{if}\ a<0,\\
\log(T+2R) & \mbox{if}\ a=0,\\
1 & \mbox{if}\ a>0.
\end{array}
\right.
\end{equation}
\end{prop}
The proof of Proposition \ref{prop:apriori} is established in the next section.
Set
\[
M:=\|f'\|_{L^\infty(\R)}+\|f''\|_{L^\infty(\R)}+\|g\|_{L^\infty(\R)}+\|g'\|_{L^\infty(\R)}.
\]

\vskip10pt
\par\noindent
{\bf The convergence of the sequence $\v{\{U_j\}}$.}
\par
First we note that $\|U_1\|\le M\e$ by (\ref{linear_derivative}).
Since (\ref{U_j}) and (\ref{apriori}) yield that
\[
\begin{array}{ll}
\|U_{j+1}\|
&\le M\e+\|L_a'(|U_j|^p)\|\\
&\le M\e+CE_a(T)\|U_j\|^p,
\end{array}
\]
the boundedness of $\{U_j\}$;
\begin{equation}
\label{bound_U}
\|U_j\|\le 2M\e\quad(j\in\N)
\end{equation}
follows from
\begin{equation}
\label{condi1}
CE_a(T)(2M\e)^p\le M\e.
\end{equation}
Assuming (\ref{condi1}), one can estimate $U_{j+1}-U_j$ as follows.
\[
\begin{array}{ll}
\|U_{j+1}-U_j\|
&\le\|L'_a(|U_j|^p-|U_{j-1}|^p)\|\\
&\le p\|L_a'\left((|U_j|^{p-1}+|U_{j-1}|^{p-1})|U_j-U_{j-1}|\right)\|\\
&\le pCE_a(T)(\|U_j\|^{p-1}+\|U_{j-1}\|^{p-1})\|U_j-U_{j-1}\|\\
&\le pCE_a(T)2(2M\e)^{p-1}\|U_j-U_{j-1}\|.
\end{array}
\]
Therefore the convergence of $\{U_j\}$ follows from
\begin{equation}
\label{convergence}
\|U_{j+1}-U_j\|\le\frac{1}{2}\|U_j-U_{j-1}\|\quad(j\ge2)
\end{equation}
provided (\ref{condi1}) and
\begin{equation}
\label{condi2}
pCE_a(T)2(2M\e)^{p-1}\le\frac{1}{2}
\end{equation}
are fulfilled.

\vskip10pt
\par\noindent
{\bf The convergence of the sequence $\v{\{(U_j)_x\}}$.}
\par
First we note that $\|(U_1)_x\|\le M\e$ by (\ref{linear_2derivative}).
Assume that (\ref{condi1}) and (\ref{condi2}) are fulfilled.
It follows from (\ref{U_j_x}) and (\ref{apriori}) that
\[
\begin{array}{ll}
\|(U_{j+1})_x\|
&\le M\e+\|L_a'\left(|U_j|^{p-1}|(U_j)_x|\right)\|+|1+a|\|L'_{a+1}(|U_j|^p)\|\\
&\le M\e+CE_a(T)\|U_j\|^{p-1}\|(U_j)_x\|+|1+a|CE_{a+1}(T)\|U_j\|^p\\
&\le M\e+CE_a(T)(2M\e)^{p-1}\|(U_j)_x\|+|1+a|CE_{a+1}(T)(2M\e)^p.
\end{array}
\]
Here we have employed the fact that
(\ref{nonlinear_derivative}) yields
\[
\begin{array}{ll}
\left|L_{a+2}'(|U_j|^px)(x,t)\right|
&\d\le\frac{1}{2}\int_0^t\frac{|U_j(x+t-s,s)|^p|x+t-s|}{\{1+(x+t-s)^2\}^{(3+a)/2}}ds\\
&\quad+\d\frac{1}{2}\int_0^t\frac{|U_j(x-t+s,s)|^p|x-t+s|}{\{1+(x-t+s)^2\}^{(3+a)/2}}ds\\
&\le L_{a+1}'(|U_j|^p)(x,t).
\end{array}
\]
Hence the boundedness of $\{(U_j)_x\}$;
\begin{equation}
\label{bound_U_x}
\|(U_j)_x\|\le 2M\e\quad(j\in\N)
\end{equation}
follows from
\begin{equation}
\label{condi3}
CE_a(T)(2M\e)^p+|1+a|CE_{a+1}(T)(2M\e)^p\le M\e.
\end{equation}
Assuming (\ref{condi3}), one can estimate $(U_{j+1})_x-(U_j)_x$ as follows.
\[
\begin{array}{ll}
\|(U_{j+1})_x-(U_j)_x\|
&\le\|L'_a(|U_j|^{p-2}U_j(U_j)_x-|U_{j-1}|^{p-2}U_{j-1}(U_{j-1})_x)\|\\
&\quad+|1+a|\|L'_{a+2}\left((|U_j|^p-|U_{j-1}|^p)x\right)\|.
\end{array}
\]
The first term on the right hand side of above inequality is split into three pieces
according to
\[
\begin{array}{l}
|U_j|^{p-2}U_j(U_j)_x-|U_{j-1}|^{p-2}U_{j-1}(U_{j-1})_x\\
=(|U_j|^{p-2}-|U_{j-1}|^{p-2})U_j(U_j)_x\\
\quad+|U_{j-1}|^{p-2}(U_j-U_{j-1})(U_j)_x\\
\quad+|U_{j-1}|^{p-2}U_{j-1}((U_j)_x-(U_{j-1})_x).
\end{array}
\]
Since
\[
\begin{array}{l}
\left||U_j|^{p-2}-|U_{j-1}|^{p-2}\right|
\\
\le
\left\{
\begin{array}{ll}
(p-2)(|U_j|^{p-3}+|U_{j-1}|^{p-3})|U_j-U_{j-1}| & \mbox{when}\ p\ge3,\\
|U_j-U_{j-1}|^{p-2} & \mbox{when}\ 2<p<3,\\
0 & \mbox{when}\ p=2,
\end{array}
\right.
\end{array}
\]
the similar manner of  handling $L_{a+2}'$ to above computations leads to
\[
\begin{array}{l}
\|(U_{j+1})_x-(U_j)_x\|\\
\le CE_a(T)\|U_j\|\|(U_j)_x\|\times\\
\qquad\times
\left\{
\begin{array}{ll}
(p-2)(\|U_j\|^{p-3}+\|U_{j-1}\|^{p-3})\|U_j-U_{j-1}\| & \mbox{when}\ p\ge3,\\
\|U_j-U_{j-1}\|^{p-2} & \mbox{when}\ 2<p<3,\\
0 & \mbox{when}\ p=2
\end{array}
\right.
\\
\quad+CE_a(T)\|U_{j-1}\|^{p-2}\|U_j-U_{j-1}\|\|(U_j)_x\|\\
\quad+CE_a(T)\|U_{j-1}\|^{p-1}\|(U_j)_x-(U_{j-1})_x\|\\
\quad+|1+a|CE_{a+1}(T)p(\|U_j\|^{p-1}+\|U_{j-1}\|^{p-1})\|U_j-U_{j-1}\|.
\end{array}
\]
Hence it follows from (\ref{convergence}) that
\[
\begin{array}{ll}
\|(U_{j+1})_x-(U_j)_x\|
&\le CE_a(T)(2M\e)^{p-1}\|(U_j)_x-(U_{j-1})_x\|\\
&\quad\d+
\left\{
\begin{array}{ll}
\d O\left(\frac{1}{2^{(p-2)j}}\right) & \mbox{when}\ 2<p<3,\\
\d O\left(\frac{1}{2^j}\right) & \mbox{otherwise}
\end{array}
\right.
\end{array}
\]
as $j\rightarrow\infty$.
Here we have employed the fact that $E_{a+1}(T)$
is dominated by $E_a(T)$ with some positive constant.
Therefore we obtain the convergence of $\{(U_j)_x\}$ provided
\begin{equation}
\label{condi4}
CE_a(T)(2M\e)^{p-1}\le\frac{1}{2}.
\end{equation}

\vskip10pt
\par\noindent
{\bf Continuation of the proof.}
\par
It is easy to find a positive constant $C_0$ independent of $\e$ and $T$ such that
all the conditions, (\ref{condi1}), (\ref{condi2}), (\ref{condi3}), (\ref{condi4}),
on the convergence of $\{U_j\}$ in the closed subspace of $X$
satisfying $\|U\|,\|U_x\|\le 2M\e$
 follows from
\[
C_0\e^{p-1}E_a(T)\le1.
\]
Therefore we obtain Theorem \ref{thm:lower-bound} for $p\ge2$.

\par
For $1<p<2$, $X$ and $M$ in the proof for $p\ge2$ above are replaced
with $Y$ and $N$ respectively, where 
\[
Y:=\{U\in C(\R\times[0,T]) :\ \mbox{supp}\ U\subset\{|x|\le t+R\}\}
\]
equipping
\[
\|U\|_Y:=\|U\|=\sup_{\R\times[0,T]}|U(x,t)|
\]
and
\[
N:=\|f'\|_{L^\infty(\R)}+\|g\|_{L^\infty(\R)},
\]
respectively.
It is trivial that the convergence of $\{U_j\}$ in the closed subspace $Y$ 
satisfying $\|U\|\le 2N\e$ follows from (\ref{condi1}) and (\ref{condi2}),
so that the proof of Theorem \ref{thm:lower-bound} is completed now
by taking $\e$ small enough.
\hfill$\Box$

%%%%%%%%%%%%%%%%%%%%%%%%%%%%%%%%%%%%%%%%%%%%%%%%%%
%%%%%%%%%%%%%%%%%%%%% SECTION4 %%%%%%%%%%%%%%%%%%%%%%%
%%%%%%%%%%%%%%%%%%%%%%%%%%%%%%%%%%%%%%%%%%%%%%%%%%

\section{Proof of Proposition \ref{prop:apriori}}
\par
In this section, we prove a priori estimate (\ref{apriori}).
Recall the definition of $L'_a$ in (\ref{nonlinear_derivative}).
From now on, a positive constant $C$ independent of $T$ and $\e$
may change from line to line.
Since
\[
\frac{1}{2}(1+|x|)\le(1+x^2)^{1/2}\le(1+|x|),
\]
we have that
\[
|L_a'(|U|^p)(x,t)|\le
C\|U\|^p\{I_+(x,t)+I_-(x,t)\},
\]
where the integrals $I_+$ and $I_-$ are defined by
\[
I_\pm(x,t):=\int_0^t\frac{\chi_\pm(x,t;s)}{(1+|t-s\pm x|)^{1+a}}ds
\]
and the characteristic functions $\chi_+$ and $\chi_-$ are defined by
\[
\begin{array}{ll}
\chi_\pm(x,t;s)
&:=\chi_{\{s: |t-s\pm x|\le s+R\}}\\
&=
\left\{
\begin{array}{ll}
1 & \mbox{when $s$ satisfies }|t-s\pm x|\le s+R,\\
0 & \mbox{otherwise},
\end{array}
\right.
\end{array}
\]
respectively.
First we note that it is sufficient to  estimate $I_\pm$ for $x\ge0$ due to its symmetry,
\[
I_+(-x,t)=I_-(x,t).
\]
Hence it follows from $0\le x\le t+R$ as well as
\[
|t-s+x|\le s+R\quad\mbox{and}\quad 0\le s\le t
\]
that
\[
\frac{t+x-R}{2}\le s\le t,
\]
so that
\[
I_+(x,t)\le\int_{(t+x-R)/2}^t\frac{1}{(1+t-s+x)^{1+a}}ds.
\]
When $a<0$, we have
\[
I_+(x,t)\le C\left(1+t+x-\frac{t+x-R}{2}\right)^{-a}\le C(T+2R)^{-a}.
\]
When $a=0$, we have
\[
I_+(x,t)\le\log\frac{1+t+x-(t+x-R)/2}{1+x}\le\log(T+2R).
\]
When $a>0$, we have
\[
I_+(x,t)\le C(1+x)^{-a}\le C.
\]
Therefore we obtain
\[
I_+\le CE_a(T)\quad\mbox{in}\ \R\times[0,T].
\]

\par
On the other hand, the estimate for $I_-$ is divided into two cases.
If $t-x\ge0$, then  $|t-s-x|\le s+R$ yields that
\[
I_-(x,t)\le\int_{(t-x-R)/2}^{t-x}\frac{1}{(1+t-s-x)^{1+a}}ds+\int_{t-x}^t\frac{1}{(1-t+s+x)^{1+a}}ds
\]
follows.
When $a<0$, we have
\[
I_-(x,t)\le C(1+t-x-(t-x-R)/2)^{-a}+C\le C(T+2R)^{-a}.
\]
When $a=0$, we have
\[
I_-(x,t)=\log(1+t-x-(t-x-R)/2)+\log(1+x)\le2\log(T+2R).
\]
When $a>0$, we have
\[
I_-(x,t)\le C.
\]
Therefore we obtain
\[
I_-\le CE_a(T)\quad\mbox{in}\ \R\times[0,T]\cap\{t-x\ge0\}.
\]
If $(-R\le)t-x\le0$, $|t-s-x|\le s+R$ yields that
\[
I_-(x,t)\le\int_0^t\frac{1}{(1+s-t+x)^{1+a}}ds.
\]
When $a<0$, we have
\[
I_-(x,t)\le C(1+x)^{-a}\le C(T+2R)^{-a}.
\]
When $a=0$, we have
\[
I_-(x,t)\le\log\frac{1+x}{1-t+x}\le\log(T+2R).
\]
When $a>0$, we have
\[
I_-(x,t)\le C(1-t+x)^{-a}\le C.
\]
Therefore we obtain
\[
I_-\le CE_a(T)\quad\mbox{in}\ \R\times[0,T]\cap\{-R\le t-x\le0\}.
\]
Summing up all the estimates for $I_+$ and $I_-$, we have
\[
|L_a'(|U|^p)|\le
C\|U\|^pE_a(T)\quad\mbox{in}\ \R\times[0,T].
\]
This completes the proof of Proposition \ref{prop:apriori}.
\hfill$\Box$.

%%%%%%%%%%%%%%%%%%%%%%%%%%%%%%%
%%%%%%%%%%%%%%%%%%%%% SECTION5 %%%%%%%%%%%%%%%%%%%%%%%
%%%%%%%%%%%%%%%%%%%%%%%%%%%%%%%%%%%%%%%%%%%%%%%%%%

\section{Proof of Theorem \ref{thm:upper-bound}}
\par
In this section, a positive constant $C$ independent of $T$ and $\e$
may change from line to line.
Let $U\in C(\R\times[0,T])$ with
\begin{equation}
\label{supp_integralsol}
\mbox{supp}\ U(x,t)\subset\{(x,t)\in\R\times[0,T]:|x|\le t+R\}
\end{equation}
 be a solution of the integral equation (\ref{integral_derivative}),
namely
\[
U=\e u_t^0+L'_a(|U|^p).
\]
Then it is easy to see by simple integration that
\[
V(x,t):=\int_0^tU(x,s)ds+\e f(x)
\]
satisfies a integral equation,
\[
V=\e u^0+L_a(|V_t|^p).
\]

\par
Set $t=x+R,\ x\ge R$. Then, inverting the order of the $(y,s)$-integral and
diminishing its domain, we have that
\[
L_a(|V_t|^p)(x,t)\ge C\int_R^x\frac{1}{(1+y)^{1+a}}dy
\int_{y-R}^{y+R}|V_t(y,s)|^pds.
\]
Hence H\"older's inequality yields that
\[
L_a(|V_t|^p)(x,t)\ge C\int_R^x\frac{1}{(1+y)^{1+a}}dy
\left|\int_{y-R}^{y+R}V_t(y,s)ds\right|^p
\]
which implies that, due to (\ref{supp_integralsol}),
\[
L_a(|V_t|^p)(x,t)\ge C\int_R^x\frac{|V(y,y+R)|^p}{(1+y)^{1+a}}dy.
\]
On the other hand, it follows from the assumption on the support of the data  that
\[
u^0(x,x+R)=\frac{1}{2}\int_{\R}g(x)dx=:G>0\quad\mbox{for}\ x\ge R.
\]
Hence $V$ satisfies
\begin{equation}
\label{V}
V(x,x+R)> G\e+C\int_R^x\frac{|V(y,y+R)|^p}{(1+y)^{1+a}}dy
\quad\mbox{for}\ x\ge R.
\end{equation}
We note that the equality in the inequality above can be removed without loss of generality
by taking slightly smaller $G$.

\par
Now we employ the comparison argument
with a solution of the related ordinary differential equation.
Let $W$ be a solution of 
\begin{equation}
\label{W}
W(x)= G\e+C\int_R^x\frac{|W(y)|^p}{(1+y)^{1+a}}dy
\quad\mbox{for}\ x\ge R.
\end{equation}
Then we have
\begin{equation}
\label{comparison}
V(x,x+R)>W(x)\quad\mbox{for}\ x\ge R.
\end{equation}
Because $V(R,2R)>W(R)$ and the continuity of $V,W$ yield that
(\ref{comparison}) holds in the neighborhood of $x=R$.
If there exists a point
\[
x_0:=\inf\{x\ge R:V(x,x+R)=W(x)\},
\]
we immediately reach to a contradiction,
\[
0=V(x_0,x_0+R)-W(x_0)
=C\int_R^{x_0}\frac{|V(y,y+R)|^p-|W(y)|^p}{(1+y)^{1+a}}dy>0.
\]
Hence (\ref{comparison}) is true and implies that
the existence time of $V(x,x+R)$ is less than the blow-up time of $W(x)$.

\par
Therefore the conclusion of Theorem \ref{thm:upper-bound} follows by solving
the initial value problem for ordinary differential equations,
\[
\left\{
\begin{array}{l}
\d W'=\frac{C|W|^p}{(1+x)^{1+a}} \quad\mbox{in}\ [R,\infty),\\
W(R)=G\e,
\end{array}
\right.
\]
which is equivalent to (\ref{W}).
In fact, the blow-up time $X$ of $W$ has to satisfy
\[
(G\e)^{1-p}
=
\left\{
\begin{array}{ll}
\d\frac{p-1}{-a}C\left\{(1+X)^{-a}-(1+R)^{-a}\right\}
& \mbox{for}\ a<0,\\
\d(p-1)C\log\frac{1+X}{1+R}
& \mbox{for}\ a=0.
\end{array}
\right.
\]
Since (\ref{comparison}) implies that the blow-up time $T$ of $W$
has to satisfy the inequality $T\ge X+R$, 
we have the blow-up condition (\ref{upper-bound}) by taking $\e$ small enough.
To see this, if $a<0$, one can estimate $X+R$ as
\[
\begin{array}{ll}
X+R
&\d=\left\{\frac{-a}{(p-1)C}(G\e)^{1-p}+(1+R)^{-a}\right\}^{1/(-a)}-1+R\\
&\d\le\left\{\frac{-a}{(p-1)C}(G\e)^{1-p}+2(1+R)^{-a}\right\}^{1/(-a)}.
\end{array}
\]
Therefore $\e_2$ in Theorem \ref{thm:upper-bound} should be defined by
\[
\frac{-a}{(p-1)C}(G\e_2)^{1-p}=2(1+R)^{-a}
\]
because it makes
\[
X+R\le\left\{\frac{2(-a)G^{1-p}}{(p-1)C}\right\}^{1/(-a)}\e^{-(p-1)/(-a)}
\quad\mbox{for}\ 0<\e\le\e_2.
\]
Similarly, if $a=0$, one can estimate $X+R$ as
\[
\begin{array}{ll}
X+R
&\d=\exp\left\{\frac{(G\e)^{1-p}}{(p-1)C}+\log(1+R)\right\}-1+R\\
&\d\le\exp\left\{\frac{(G\e)^{1-p}}{(p-1)C}+2\log(1+R)\right\}.
\end{array}
\]
Therefore $\e_2$ in Theorem \ref{thm:upper-bound} should be defined by
\[
\frac{(G\e_2)^{1-p}}{(p-1)C}=2\log(1+R)
\]
because it makes
\[
X+R\le\exp\left\{\frac{2G^{1-p}}{(p-1)C}\e^{-(p-1)}\right\}
\quad\mbox{for}\ 0<\e\le\e_2.
\]
The proof of Theorem \ref{thm:upper-bound} is now completed.
\hfill$\Box$

%%%%%%%%%%%%%i%%%%%%%%%%%%%%%%%%%%%%%
%%%%%%%%%%%% Acknowledgement %%%%%%%%%%%%%%%
%%%%%%%%%%%%%%%%%%%%%%%%%%%%%%%%%%%%%
\section*{Acknowledgement}
\par
This work is inspired by personal communications with Professor Kunio Hidano
(Mie Univ., Japan).
All the authors appreciate him for fruitful discussions.
They also appreciate the reviewers for their many helpful comments
which make the manuscript completed.
The third author is partially supported
by the Grant-in-Aid for Scientific Research (A) (No.22H00097) and (B) (No.18H01132), 
Japan Society for the Promotion of Science.

%%%%%%%%%%%%%%%%%%%%%%%%%%%%%%%%%%%%%%
%%%%%%%%%%%% References %%%%%%%%%%%%%%%%%%%%
%%%%%%%%%%%%%%%%%%%%%%%%%%%%%%%%%%%%%%

\bibliographystyle{plain}

\begin{thebibliography}{20}

\bibitem{John79}{F. John},
{\it Blow-up of solutions of nonlinear wave equations in three space dimensions},
Manuscripta Math., {\bf 28} (1979), 235-268.

\bibitem{John_book}{F. John},
\lq\lq Nonlinear Wave Equations, Formation of Singularities",
ULS Pitcher Lectures in Mathematical Science, Lehigh University,
American Mathematical Society, Providence, RI, 1990.

\bibitem{KTW19} M. Kato, H. Takamura and K. Wakasa,
{\it The lifespan of solutions of semilinear wave equations
with the scale-invariant damping in one space dimension},
Differential Integral Equations, {\bf 32} (2019), no. 11-12, 659-678.

\bibitem{KMT21}{S. Kitamura, K. Morisawa and H. Takamura},
{\it The lifespan of classical solutions of semilinear wave equations
with spatial weights and compactly supported data
in one space dimension},
J. Differential Equations, {\bf 307} (2022), 486-516.

\bibitem{KOY13}{H. Kubo, A. Osaka and M. Yazici},
{\it Global existence and blow-up for wave equations
with weighted nonlinear terms in one space dimension},
Interdisciplinary Information Sciences, {\bf 19} (2013), 143-148.

\bibitem{LYZ91}
T.-T. Li (D.-Q. Li), X. Yu and Y. Zhou,
{\it Dur\'ee de vie des solutions r\'eguli\`eres pour les \'equations des ondes non lin\'eaires unidimensionnelles} (French),
C. R. Acad. Sci. Paris S\'er. I Math., {\bf 312} (1991), no. 1, 103-105.

\bibitem{LYZ92}
T.-T. Li, X. Yu and Y. Zhou,
{\it Life-span of classical solutions to one-dimensional nonlinear wave equations},
Chinese Ann. Math., Ser. B, {\bf13} (1992), no. 3, 266-279. 

\bibitem{Suzuki10}{A. Suzuki},
\lq\lq Global Existence and Blow-Up of solutions
to Nonlinear Wave Equations in One Space Dimension" (Japanese),
Master Thesis, Saitama University, 2010.

\bibitem{Wakasa17}{K. Wakasa},
{\it The lifespan of solutions to wave equations with weighted nonlinear terms in one space dimension},
Hokkaido Math. J., {\bf 46} (2017), 257-276.

\bibitem{Zhou92}{Y. Zhou},
{\it Life span of classical solutions to $u_{tt}-u_{xx}=|u|^{1+\alpha}$},
Chin. Ann. Math. Ser.B, {\bf 13} (1992), 230-243.

\bibitem{Zhou01}{Y. Zhou},
{\it Blow up of solutions to the Cauchy problem for nonlinear wave equations},
Chinese Ann. Math. Ser. B, {\bf22} (2001), no. 3, 275-280.

\end{thebibliography}

\end{document}